\magnification=\magstephalf
\def\hexnumber#1{\ifcase#1 0\or 1\or 2\or 3\or 4\or 5\or 6\or 7\or 8\or
 9\or A\or B\or C\or D\or E\or F\fi}
%
%
\font\twelvemsa=msam10 scaled 1200   
\font\tenmsa=msam10                  
\font\ninemsa=msam9            \font\sevenmsa=msam7
\font\sixmsa=msam6             \font\fivemsa=msam5
%
%
\newfam\msafam                 \textfont\msafam=\tenmsa
\scriptfont\msafam=\sevenmsa   \scriptscriptfont\msafam=\fivemsa
\edef\hexa{\hexnumber\msafam}        
\def\msa{\fam\msafam\tenmsa}         
%
%
\font\twelvemsb=msbm10 scaled 1200   
\font\tenmsb=msbm10                  
\font\ninemsb=msbm9            \font\sevenmsb=msbm7
\font\sixmsb=msbm6             \font\fivemsb=msbm5
%
\newfam\msbfam                 \textfont\msbfam=\tenmsb
\scriptfont\msbfam=\sevenmsb   \scriptscriptfont\msbfam=\fivemsb
\edef\hexb{\hexnumber\msbfam}        
\def\msb{\fam\msbfam\tenmsb}         
%
%
\font\twelveeufm=eufm10 scaled 1200  
\font\teneufm=eufm10                 
\font\nineeufm=eufm9           \font\seveneufm=eufm7
\font\sixeufm=eufm6            \font\fiveeufm=eufm5
%
\newfam\eufmfam                \textfont\eufmfam=\teneufm
\scriptfont\eufmfam=\seveneufm \scriptscriptfont\eufmfam=\fiveeufm
\edef\hexf{\hexnumber\eufmfam}      
\def\frak{\fam\eufmfam\teneufm}     
%
%
%
\font\twelverm=cmr10 scaled 1200    
\font\ninerm=cmr9                   
\font\sixrm=cmr6
%
\font\twelvei=cmmi10 scaled 1200    
\font\ninei=cmmi9                   
\font\sixi=cmmi6
%
\font\twelvesy=cmsy10 scaled 1200   
\font\ninesy=cmsy9                  
\font\sixsy=cmsy6
%
\font\twelvebf=cmbx10 scaled 1200   
\font\ninebf=cmbx9                  
\font\sixbf=cmbx6
%
%
\font\twelveit=cmti10 scaled 1200   
\font\nineit=cmti9                  
%
\font\twelvesl=cmsl10 scaled 1200   
\font\ninesl=cmsl9                  
%
\font\twelvett=cmtt10 scaled 1200   
\font\ninett=cmtt9                  
%
%
%
%
\def\small{%
%
%
\textfont0=\ninerm \scriptfont0=\sixrm \scriptscriptfont0=\fiverm
\def\rm{\fam0\ninerm}        
%
%
\textfont1=\ninei \scriptfont1=\sixi \scriptscriptfont1=\fivei
%
%
\textfont2=\ninesy \scriptfont2=\sixsy \scriptscriptfont2=\fivesy
%
%
\textfont3=\tenex \scriptfont3=\tenex \scriptscriptfont3=\tenex
%
%
\textfont\bffam=\ninebf \scriptfont\bffam=\sixbf
\scriptscriptfont\bffam=\fivebf \def\bf{\fam\bffam\ninebf}%
%
%
\textfont\itfam=\nineit \def\it{\fam\itfam\nineit}%
\textfont\slfam=\ninesl \def\sl{\fam\slfam\ninesl}%
\textfont\ttfam=\ninett \def\tt{\fam\ttfam\ninett}%
%
%
%
\textfont\msafam=\ninemsa \scriptfont\msafam=\sixmsa
\scriptscriptfont\msafam=\fivemsa \def\msa{\fam\msafam\ninemsa}%
%
%
\textfont\msbfam=\ninemsb \scriptfont\msbfam=\sixmsb
\scriptscriptfont\msbfam=\fivemsb \def\msb{\fam\msbfam\ninemsb}%
%
%
\textfont\eufmfam=\nineeufm  \scriptfont\eufmfam=\sixeufm
\scriptscriptfont\eufmfam=\fiveeufm \def\frak{\fam\eufmfam\nineeufm}%
%
%
%
\normalbaselineskip=11pt
\setbox\strutbox=\hbox{\vrule height8pt depth3pt width0pt}%
%
%
\normalbaselines\rm}    
%
%
%
%
\def\large{%
\textfont0=\twelverm \scriptfont0=\ninerm \scriptscriptfont0=\sevenrm
\def\rm{\fam0\twelverm}%
\textfont1=\twelvei \scriptfont1=\ninei \scriptscriptfont1=\seveni
\textfont2=\twelvesy \scriptfont2=\ninesy \scriptscriptfont2=\sevensy
\textfont3=\tenex \scriptfont3=\tenex \scriptscriptfont3=\tenex
\textfont\bffam=\twelvebf \scriptfont\bffam=\ninebf
\scriptscriptfont\bffam=\sevenbf \def\bf{\fam\bffam\twelvebf}%
\textfont\itfam=\twelveit \def\it{\fam\itfam\twelveit}%
\textfont\slfam=\twelvesl \def\sl{\fam\slfam\twelvesl}%
\textfont\ttfam=\twelvett \def\tt{\fam\ttfam\twelvett}%
\textfont\msafam=\twelvemsa \scriptfont\msafam=\ninemsa
\scriptscriptfont\msafam=\sevenmsa \def\msa{\fam\msafam\twelvemsa}
\textfont\msbfam=\twelvemsb \scriptfont\msbfam=\ninemsb
\scriptscriptfont\msbfam=\sevenmsb \def\msb{\fam\msbfam\twelvemsb}
\textfont\eufmfam=\twelveeufm  \scriptfont\eufmfam=\nineeufm
\scriptscriptfont\eufmfam=\seveneufm \def\frak{\fam\eufmfam\teneufm}
\normalbaselineskip=15pt
\setbox\strutbox=\hbox{\vrule height11pt depth4pt width0pt}%
\normalbaselines\rm}%
%
\def\Bbb{\msb}

%
\def\N{{\frak N}}

%
\mathchardef\plussquare="0\hexa01
\mathchardef\nge="3\hexb0B
\mathchardef\maltesecross="0\hexa7A
\mathchardef\del="0\hexf01
%

%

\input epsf
\input colordvi

\def\to{\rightarrow\ }

\def\nl{\hfill\break}
\def\s{\sigma}

\def\smaths#1{$\textfont0=\sevenrm \textfont1=\seveni \scriptfont0=\fiverm \scriptfont1=\fivei #1$}

\def\pagebreak{\vfill\eject}

\input colordvi

\overfullrule=0pt

\font\npt=cmr9
\font\Bbb=msbm10
\font\chapfont=cmbx12
\font\secfont=cmss12 

\font\nam=cmr8
\font\aff=cmti8


\mathchardef\square="0\hexa03
\def\qed{\hfill$\square$\par\rm}
\def\boxing#1{\ \lower 3.5pt\vbox{\vskip 3.5pt\hrule \hbox{\strut\vrule
\ #1 \vrule} \hrule} }
\def\up#1{ \ \vbox{\vskip 3.5pt\hrule \hbox{\strut \  #1 \vrule} } }

\def\down#1{\ \lower 3.5pt\vbox{\vskip 3.5pt \hbox{\strut \ #1 \vrule} \hrule} }
\def\negdown#1{\ \lower 3.5pt\vbox{\vskip 3.5pt \hbox{\strut  \vrule \ #1 }\hrule} }

\def\items{\par\leftskip = 25pt \parskip=0.5\baselineskip}
\def\enditems{\par\leftskip = 0pt \parskip=\baselineskip}

\newcount\diagramnumber \diagramnumber=0
\newcount\chapternumber \chapternumber=1
\newcount\questionnumber \questionnumber=0  
\newcount\secnum \secnum=0
\newcount\subsecnum \subsecnum=0
\newcount\subsubsecnum \subsubsecnum=0
\newcount\defnum \defnum=0 
\newcount\theonum \theonum=0
\newcount\lemnum \lemnum=0
\newcount\cornum \cornum=0
\newcount\exnum \exnum=0
\newcount\figurenumber \figurenumber=0
\newcount\bibnum \bibnum=0

\def\chapter#1{\global\advance\chapternumber by 1
\np\diagramnumber=0\questionnumber=0
\centerline{\Large\noindent
{\number\chapternumber\ \ }#1}         
\ppar                                  
\resultnumber=1                        
\def\chaptername{#1}\tenrm}       

\hsize=5.5 truein
\vsize=9.25 truein
\voffset=0.5 in

\baselineskip=.5 truecm
\parskip=\baselineskip
 1

\parindent=0pt

\def\Z{\hbox{\Bbb Z}}
\def\R{\hbox{\Bbb R}}

\def\N{\hbox{\Bbb N}}

\def\x{\hbox{\bf x}}

\def\oover#1{\vbox{\ialign{##\crcr
{\npt o}\crcr\noalign{\kern 1pt\nointerlineskip}
$\hfil\displaystyle{#1}\hfil$\crcr}}}

\newif \iftitlepage \titlepagetrue


\def\diagram#1{\global\advance\diagramnumber by 1
$$\epsfbox{birack-polynomialfig.\number\diagramnumber}$$\postdisplaypenalty=10000
\centerline{\npt Fig.\ \the\diagramnumber. \npt\sl\ #1}\par}

\def\section#1{
                \vskip 10 pt
                \advance\secnum by 1 \subsecnum=0 \subsubsecnum=0
                \leftline{\secfont \the\secnum \rm\quad\secfont #1}
                }

\def\subsection#1{
                \vskip 10 pt
                \advance\subsecnum by 1 \subsubsecnum=0
                 \leftline{\secfont \the\secnum.\the\subsecnum\ \rm\quad\secfont \ #1}
                }

\def\subsubsection#1{
                \vskip 10 pt
                \advance\subsubsecnum by 1 
                 \leftline{\secfont \the\secnum.\the\subsecnum.\the\subsubsecnum \rm\quad\secfont \ #1}
                }

 \def\theorem#1{
                \advance\theonum by 1 
                \par\bf Theorem  \the\secnum.\the\theonum \sl\ #1\rm\nl
               }

\def\lemma#1{
                \advance\lemnum by 1 
                \par\bf Lemma  \the\secnum.\the\lemnum \sl\ #1\rm \nl 
                }

\def\corollary#1{
                \advance\cornum by 1 
                \par\bf Corollary  \the\secnum.\the\cornum \sl\ #1\rm\ \par
               }

\def\example#1{
                \advance\exnum by 1 
                \par\bf Example  \the\secnum.\the\exnum \sl\ #1  \par\rm
               }
               
\def\definition#1{
                \advance\defnum by 1 
                \par \bf Definition \the\secnum .\the\defnum  \it \ #1\rm \par
                }

\def\cite#1{{\bf [#1]}}\rm

\def\chaptername{  }

\def\today{\number\day\space\ifcase\month\or January\or February\or
March\or April\or May\or June\or July\or August\or September\or
October\or November\or December\fi\space\number\year}


\today
\vglue 20 pt
\centerline{\bf Generalised Biracks and the Birack Polynomial Invariant}
\medskip
\centerline{\chapfont Dedicated to the memory of Roger Fenn}
\medskip
\vglue 20 pt
\centerline{Andrew Bartholomew${}^1$, Roger Fenn${}^2$ 
Louis Kauffman${}^3$,} 
\centerline{\aff ${}^1$, ${}^2$School of Mathematical Sciences, University of Sussex}
\centerline{\aff Falmer, Brighton, BN1 9RH, England}

\centerline{\aff ${}^3$Department of Mathematics, Statistics and 
Computer Science, 851 South Morgan St}
\centerline{\aff University of Illinois at Chicago}
\centerline{\aff Chicago, Illinois, USA 60607-7045}

\bigskip

\medskip
%
\centerline{\nam ABSTRACT}
\leftskip=0.5 in


{\bf  We generalise the finite biquandle colouring invariant to a polynomial invariant based on labelling a knot diagram with a finite birack that reduces to the biquandle 
colouring invariant in that case.  The 
polynomial is an invariant of a class of knot theories amenable to a generalisation of theorem of Trace on regular homotopy.  We take the 
opportunity to reprise the relevant generalised knot theory and the theory of generalised biracks in the light of this work and recent developments.}



\hsize=6.125 truein
\vsize=9.25 truein
\voffset=0.5 in
\leftskip=0in
\baselineskip=.5 truecm
\parskip=\baselineskip

\section{Introduction}

From the outset (\cite{FJK}) the theory of biracks and biquandles has been applied to both classical and virtual knots and links and they 
were given their first systematic treatment in a generalised knot theory in \cite{F1}.  Here we build upon this work, taking advantage of subsequent insights, to present
an account of an invariant of a set of knot theories that is based on labelling a diagram with a birack, without the need for
separate names and definitions for each theory.

The usual fundamental expression of a knot is a representative diagram within an equivalence class. Here we take 
as the fundamental expression of a knot $K$ the set of {\it all} diagrams representing $K$, which is written  as ${\cal D}(K)$ and can be
thought of as an infinite graph with vertices the diagrams and edges labelled by the Reidemeister moves. Alternatively, it may be viewed as a groupoid category,
with objects representative diagrams and morphisms generated by Reidemeister moves.

In the case of classical knots, ${\cal D}(K)$ has a 2-fold subdivision into subsets ${\cal D}_{w,n}(K)$.
The vertices of ${\cal D}_{w,n}(K)$ are the diagrams  with {\it writhe} $w$ and {\it turning number} $n$ (defined below).  The edges of ${\cal D}_{w,n}(K)$ are labelled by one of the second or third Reidemeister moves. In other words the diagrams of ${\cal D}_{w,n}(K)$ are connected by a {\it regular homotopy}.
By a theorem of Trace, \cite{T}, a regular homotopy class is determined by the pair $w,n$. It follows that the ${\cal D}_{w,n}(K)$ are connected subgraphs of ${\cal D}(K)$ and the neighbours are joined by edges labelled by the first Reidemeister move.

By generalising Trace's result we are able to extend this subdivision to a number of different knot theories.  In the classical case, 
a consequence of the bi-division of ${\cal D}(K)$ is that any invariant of regular homotopy is the same for any representative diagram in  
${\cal D}_{w,n}(K)$ and we have the corresponding property in our more general setting.

The contents of this paper are as follows. In section 2 we provide a re-evaluation of generalised knots, and in sections 3 and 4, of biracks 
and biquandles.  In section 5 we generalise Trace's theorem of classical knots and discuss the implications for the graph ${\cal D}(K)$.
In section 6 we define our invariant, called the {\it birack polynomial} and give concrete examples of the polynomial for some simple
classical knots, showing that it allows us to distinguish the trefoil and figure eight knots and can prove that they are non-trivial.  We also 
show that the birack polynomial is able to distinguish as non-trivial the two Bigelow braids used in \cite{B} to demonstrate that the 
Burau representation of the braid group $B_n$ is not failthful for $n=5,6$.

\section{Generalised knot theories and $R$ moves}

The $R$ stands for Reidemeister. His moves are generalised for a generalised knot theory. Full details can be found in \cite{F1} however for this paper we just need a brief description as follows. 

Objects in a generalised knot theory, $\cal K$, can be defined combinatorially in terms of a diagram, $D$ involving one or more
crossing types $C_t$.  Each crossing type has a {\it positive} and {\it negative} version, called the {\it polarity} which may not be 
distinct.  We represent crossing types with letters, using $r, v, f, s$ for the familiar real, virtual, flat and singular crossings respectively 
and letters $a,b,\ldots$ for general crossing types.  We use a bar, $\bar a$, to represent a crossing of type $a$ having the opposite polarity to $a$ 
and say that a crossing labelled with crossing type $a$ is {\it tagged} by $a$ and that $a$ is a crossing {\it tag}.

Generalised theories such as multi-virtual knot theory, described in \cite{K2}, where there is a multiplicity of virtual crossing types are accommodated by
regarding each type of virtual crossing as a crossing type in its own right.

The underlying graph, $|D|$, of a diagram is the image of one or more immersed oriented circles in general position in the 
plane, so is a four-valent planar graph with additional structure.  
The edges of $|D|$ are called the {\it arcs} of the diagram.  Each vertex of $|D|$ is tagged by a signed crossing type 
included in the theory.  The objects of $\cal K$ are then equivalence classes of diagrams under an equivalence relation generated by
a collection of permitted generalised $R$ moves. What is not allowed is just as important as what is allowed.

The following $R$ moves may be specified by the theory $\cal K$.

\diagram{$R_1(a)$}
\diagram{$R_2(a)$ parallel}
\diagram{$R_2(a)$ non-parallel}
\diagram{$R_3(a,a,b)$}
\diagram{$R_4(a,b)$ commuting}

Fig.\ 1 represents four distinct $R_1$ configurations depending on orientation and the polarity of the crossing $a$.  If $\cal K$ permits
$R_1$ moves for a given crossing type it usually permits all four.  Similarly, if one is not permitted by $\cal K$, usually no variant is
and so we do not focus on the distinction between the four cases here.

If the theory permits $R_2$ moves involving an inverse pair of crossings of the type $C_t$, then the crossing type $C_t$ is called {\it 
regular} in $\cal K$.  As with $R_1$ moves, if $\cal K$ permits $R_2(a)$ it usually permits both the parallel and non-parallel versions and 
we shall assume this 
to be the case for all regular crossing types.  There is an additional non-parallel configuration to that shown in Fig.\ 3 where the orientation of the arcs is 
reversed and the bigon oriented anti-clockwise rather than clockwise, which we also assume to be permitted if $\cal K$ permits $R_2(a)$.  As with $R_1(a)$, if $\cal K$
does not permit $R_2(a)$ is usually does not permit any of these variants.

A knot theory containing a regular crossing type is itself called regular.  We shall only consider regular knot theories in what follows.

It is shown in \cite{F1} that we need consider only $R_3$ moves where the three strands are oriented braid-like as in Fig.\ 4. As indicated by the 
notation, the move describes how two signed crossings of the same type and polarity may be considered to move past a signed crossing of a different type.
If the move $R_3(a,a,b)$ is permitted, we say that $a$ {\it dominates} $b$.

The $R_4(a,b)$ move is permitted in singular knot theories and is the basis of the Kauffman virtualisation move between a 
virtual crossing and either a classical or flat crossing.  Like $R_2$, it comes in a parallel and non-parallel version.

\subsection {Example knot theories}

Classical knot theory has a single crossing type, $r$ and permits $R_1(r), R_2(r), R_3(r,r,r)$.

Virtual knot theory has crossing types $r , v$ and permits $R_1(r), R_2(r), R_3(r,r,r)$, \break 
$R_1(v), R_2(v), R_3(v,v,v)$ and $R_3(v,v,r)$.

Flat virtual knot theory has crossing types $f$ and $v$ and permits the same $R$ moves as virtual knot theory but with the $f$ tag replacing the $r$ tag.

Each of the above theories $\cal K$, has a singular variety $\cal SK$ by adding the singular crossing type $s$ and permitting $R_2(s,s)$, 
$R_3(a,a,s)$, $R_3(\bar a, \bar a,s)$ and $R_4(a,s)$ for each of the crossing types $a \ne s$ in $\cal K$.

\section{Square functions, $X^2\to X^2$.} 

Let $X$ be a set, called the {\it label set}, and $F:X^2 \rightarrow X^2$ be a map, which we assume to respect any structure on $X$. A useful notation for $F$ is 
$$F(x,y)=(f_x(y), f^y(x))$$ 
where $f_x:X\to X$ is called the {\it down} map indexed by $x$ and $f^y:X\to X$ is called the {\it up} map indexed by $y$. The up and down maps are called the {\it auxiliary} functions of $F$.

We can illustrate the situation by placing the variables at the corner of a square  as in Fig.\ 6, where we show the two examples,
$S$ and $T$, of $F$ that are of particular interest. 

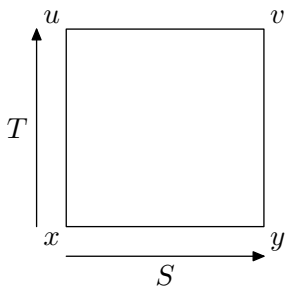
\diagram{The mapping square}

In order to derive invariants, the arcs of knot diagram are labelled with elements of $X$ in such a way that the relationship between the 
four labels at a crossing is described by either $S$ or $T$.  Traditionally, the focus has been on the map $S$ but there are situations 
where considering $T$ is useful and, as we show below, in most knot theories, the two views are equivalent.

We call a map $S$, taking the column ${u\choose x}$ to the column ${v\choose y}$ a {\it switch} and a  map $T$, taking the row 
$(x,y)$ to the row  $(u,v)$ a {\it twitch}. We shall continue to use the the generic notation $F$ to refer to either $S$ or $T$ when we 
wish not to be specific.

%

A function, $F:X^2\to X^2$, is called {\it formed} if the auxiliary functions  are always bijections, and {\it fully formed} if in addition $F$ is also a bijection. It is easy to find examples of formed functions which are not fully formed.

\lemma{ \nl
1. If $T$ ($S$) is formed then $S$ ($T$) is well defined and is a bijection.\nl
2. $T$ is fully formed if, and only if, $S$ is fully formed.}
{\bf Proof.}\nl
1. Using the notation above, if $T(x,y)$ is formed then
$$S{u\choose x}=S{t_x(y)\choose x}={t^y(x)\choose y}={t^{t_x^{-1}(u)}(x)\choose t_x^{-1}(u)}$$
This is a bijection because, given any $(t^y(x),y)$, $x$ is uniquely defined and therefore so is $t_x(y)$. So the inverse given by 
$$S^{-1}{t^y(x)\choose y} = {t_x(y)\choose x}$$
 is uniquely defined.
A similar argument works if $S$ is formed.\nl
2. From statement 1, for a given $T$ it remains to show that $S$ is formed.  Using the functional notation for S we have $s_x(t_x(y)) = y$, so $s_x = t_x^{-1}$ and so $s_x$ is a bijection.
For $s^y$ we proceed by contradiction.  Suppose first, that for some $a$, $s^a$ is not surjective.  Then there is a $w$ for which there is no $(a,b)$ with $s^a(b) = w$.  
Since $T$ is a bijection there is a unique pair $(x,y)$ such that $T(x,y)=(a,w)$, whereupon $S(x,a)=(y,w)$, so $s^a(x)=w$ giving a contradiction.

Now suppose $s^a$ is not injective.  Then there exist $x,x' \in X$ such that $s^a(x)=s^a(x') = b$ say. We therefore have $S(x,a) = (s_x(a),b)$
and $S(x',a) = (s_{x'}(a),b)$, hence $T(x,s_x(a)) = (a,b)$ and $T(x',s_{x'}(a)) = (a,b)$, so 
$(x,s_x(a)) = (x',s_{x'}(a))$, since $T$ is a bijection.  Thus $x=x'$ and again we have a contradiction. 
A similar argument works for a given $S$.\qed

As noted above, if we assume that the the up and down operations are bijections it does not mean that $S$ or $T$ are invertible. 
The following is an example of such a twitch $T$, and shows that the switch $S$ it defines is invertible but not formed.

\example{Let $X=\R$ and $t^y(x)=t_y(x)=x-y$. So $T(x,y)=(u,v)=(t_x(y),t^y(x))=(y-x,x-y)$ is formed but not invertible.
From the mapping square in Fig.\ 5 we see that $S(x,y-x) = (y,x-y)$.  Put $a=x, b=y-x$, then $y=x+b=a+b$ and $x-y = -b$ 
and we have 
 $$S{b\choose a}={-b\choose a+b}$$
So $S$ is invertible but not formed.}
 
A similar example shows a similar result with $S, T$ interchanged.

\section{Labelling diagrams: biracks and biquandles}

We wish to label the arcs of a diagram with elements of $X$ in such a way that the relationship between the four labels at a crossing
is described by either $S$ or $T$.

To do so we must consider the asymmetry of the mapping square introduced by the positive and negative form of a crossing type.  In Fig.\ 7 we show a positive real crossing on the left square 
and a negative real crossing on the right square in Fig.\ 7. A similar situation applies to any type of crossing in a general knot theory, but we show oriented real crossings as they have the advantage of having distinct positive and negative versions.
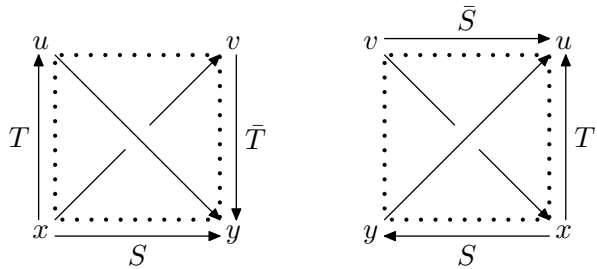
\diagram{Neighbourhoods of crossings}

The four labels, $x,y,u,v$, at a positive crossing of some type are described by the relationship $S(x,u)=(y,v)$ or $T(x,y)=(u,v)$.  We may 
describe the labels at a negative crossing of that type in the same way if the labels associated with adjacent ingress and egress arcs at
the crossing are interchanged, as shown on the right of Fig.\ 7.   

From the perspective of the mapping $S$, it is natural and convenient to consider the inverse map, denoted $\bar S$ and shown in Fig.\ 7, 
which is uniquely defined if $S$ is a bijection.  Similarly we consider $\bar T$, shown on the left of 
Fig.\ 7 to emphasise that it indicates the map $\bar T(u,v) = (x,y)$.

When a general knot theory $\cal K$, has more than one crossing type, the maps $S$ and $T$ are comprised of a set of maps, one for each crossing type 
included in the theory, so that, writing $F$ for either $S$ or $T$, 
$F = \{F_t : C_t \in {\cal K}\}$.  For each crossing of type $C_t$ the corresponding component map $F_t$ describes the relationship between the labels at that crossing.

If we wish to be explicit, we describe the binary operations defined by the component map $F_t: X^2 \rightarrow X^2$  as $F_t(a,b) = (f^t_a(b), f^b_t(a))$.  If the crossing 
type is understood, we drop the $t$ from the notation as we have been doing above and write $F(a,b) = (f_a(b), f^b(a))$.
A useful extension of  the superscript notation of quandles is to add a subscript notation, so we may write $x^a=f^a(x)$ and $x_a = f_a(x)$.
As before, we think of $a$ as {\it acting } on $x$ but now also acting up or down accordingly.

If $S$ or $T$ are bijections then as already noted we denote their inverse mappings $\bar S$ and $\bar T$.  We write $(f^a)^{-1}(x) = x^{\bar a}$ and 
$(f_a)^{-1}(x) = x_{\bar a}$, so $x^{a\bar a} = x^{\bar a a} = x_{a\bar a} = x_{\bar a a} = x$.


\subsection{The implications of the $R$ moves}

The permitted Reidemeister or $R$ moves will place conditions on the transformations $S$ and $T$ for the appropriate crossing types, which we 
now explore.

We wish to label oriented arcs of a diagram with elements of $X$ as in Fig.\ 8.
\diagram{Labelled and oriented arc}
\noindent These can be related to square crossing neighbourhoods provided the orientation is coherent and the label mappings determined by $S$ or $T$ are consistent.

\subsubsection{$R_1$}

In the case of an $R_1$ move, consider Fig.\ 9 where an arc joining the top two creates a curl based at some crossing.

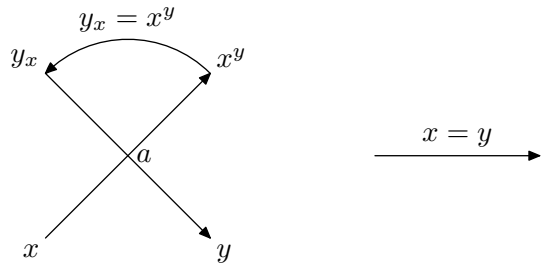
\diagram{Creating an $R_1$ curl and its deletion}

It is clear that for consistency of the labelling we must have $y_x=x^y$.

It is also clear that, if $R_1(a)$ is permitted, deleting this curl yields the right diagram of Fig.\ 9, requiring that $x=y$.
We therefore have the following:

\theorem{In order for the transformation $T_a$ to be compatible with the first Reidemeister move at a signed crossing type $a$ we need
$$x^x=x_x \hbox{ for all }x\in X \hbox{ or equivalently }T_a(\Delta)=\Delta$$
where $\Delta$ is the diagonal of $X^2$.\qed}

Equivalently, if we consider a crossing labelled by the map $S$, we have:

\theorem{In order for the transformation $S_a$ to be compatible with the first Reidemeister move at a signed crossing type $a$ we need
$$x^y=y \Leftrightarrow y_x = x \hbox{ for all }x,y\in X$$.\qed}

Note that all four of the $R_1$ variants mentioned in section 2 result in the same condition in the statement of these theorems, 
which we call the {\it biquandle condition} for the crossing type $a$.

\subsubsection{$R_2$}

Recall that there are two types of Reidemeister two moves: one where the arcs are parallel and one where they are not. Referring to Fig.\ 10 we can visualise the first by translating the right hand square so that the upright edges coincide and the squares cancel.
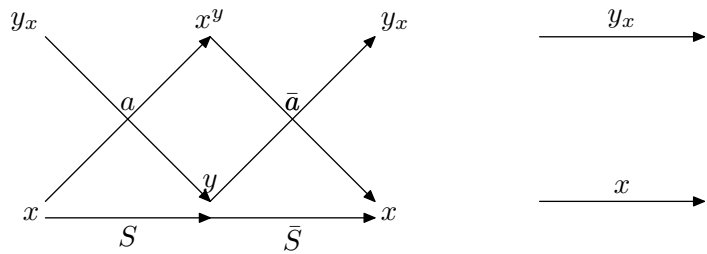
\diagram{parallel $R_2$, $S\bar S=1$}

Here the tag $a$ is a positive crossing, with $\bar a$ its negative counterpart and the crossings again labelled by the map $T$.  The position 
of the two mapping squares could be interchanged with a similar cancelling effect.

For the second kind of Reidemeister 2 move, refer to Fig.\ 11. Rotate the $\bar a$ square through $180$ degrees and place it either on top of the $a$ square or underneath.
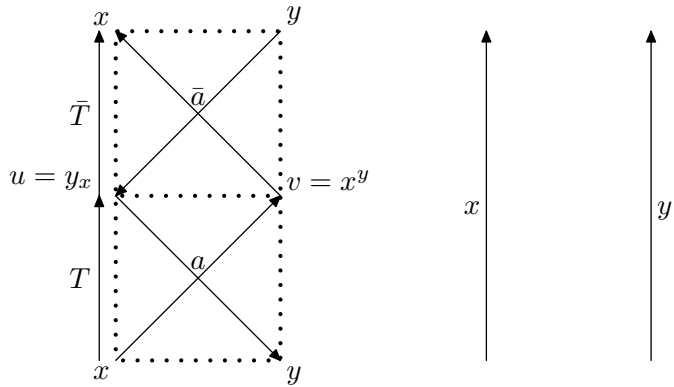
\diagram{non-parallel $R_2$, $T\bar T=1$}

In view of the transition between left and right of Fig.\ 10, and between the top and bottom of Fig.\ 11, it can be seen that these diagrams are
consistently labelled with respect to the maps $S$ and $T$ if both are bijections.  In view of Lemma 3.1 (2), the following should now be clear.

\theorem{In order to be compatible with the two kinds of Reidemeister 2 moves for crossing type $a$ we need $S_a$, or equivalently, $T_a$ to be fully formed.\qed} 

In a general knot theory, $\cal K$, each component map may independently be formed or fully formed, though in many cases, all component maps are fully formed.
Recall that we call crossing type $C_t$ {\it regular} in $\cal K$ if the theory permits $R_2$ moves involving an inverse pair of crossings of the type $C_t$ but we do not assume that all crossing types in a theory are regular.

\subsubsection{$R_3$}

Consider Fig.\ 12 where square neighbourhoods of the crossings involved in an $R_3(a,a,b)$ move are shown.

\diagram{$R_3(a,a,b)$}
If the move $R_3(a,a,b)$ is permitted, crossing type $a$ dominates crossing type $b$.  If we are to label the move with elements of $X$
so that the crossings are consistently described by the transformations $T$ and $S$ this imposes the following conditions.

\theorem{If $T_a, S_a, T_b,S_b$ are the transformations corresponding to the crossings tagged $a, b$ then in order to
be compatible with the Reidemeister 3 move  $R_3(a,a,b)$ we must have
$$\pmatrix{T_a\cr 1\cr}\pmatrix{1\cr T_a\cr}\pmatrix{T_b\cr 1\cr}=\pmatrix{1\cr T_b\cr}\pmatrix{T_a\cr 1\cr}\pmatrix{1\cr T_a\cr}$$
and also \hfill \break
\strut\hfill $(S_a\times 1)(1\times S_a)(S_b\times 1)=(1\times S_b)(S_a\times 1)(1\times S_a).$\hfil\qed}

In most knot theories each crossing type dominates itself, so the $R_3$ move $R_3(a,a,a)$ is permitted.  If we consider how the equation 
$$\pmatrix{T\cr 1\cr}\pmatrix{1\cr T\cr}\pmatrix{T\cr 1\cr}=\pmatrix{1\cr T\cr}\pmatrix{T\cr 1\cr}\pmatrix{1\cr T\cr}$$
acts on  $(x,y, z)$ then the following identities, called the $T$ {\it Yang-Baxter} identities  fall out.
$$\eqalign{
z_{xy_x}&=z_{yx^y}\cr
x^{zy^z}&=x^{yz_y}\cr
{y_x}^{z_x}&={y^z}_{x^z}\cr
}
$$
In \cite{F1} it is shown how these formul\ae\ may be summed up by putting them on the edges of a cube with a natural relationship with 
an homology theory but we do not need to consider this here. 

If we consider the move $R_3(a,a,a)$ from the perspective of $S$ we obtain the more complicated $S$ {\it Yang-Baxter} identities.
$$\eqalign{
x^{yz} &= x^{z_y y^z}\cr
x_{yz} &= x_{z^y y_z}\cr
{x_y}^{z_{y^x}} &= {x^z}_{y^{z_x}}\cr
}
$$

The simplicity of the $T$ Yang-Baxter identities is a motivation for using that perspective.


\subsubsection{$R_4$}

It is clear from Fig.\ 5 that the parallel $R_4$ moves requires that $S_a$ and $S_b$ commute.

\diagram{$R_4(a,b)$}

In the non-parallel cases we require that $T_a \bar T_b = T_b \bar T_a$ and $\bar T_a T_b = \bar T_b T_a$, see Fig.\ 13.

\subsection{Birack and biquandle definitions}

If $\cal K$ permits the dominance relation $R_3(a,a,b)$ and the map $F$ satisfies  
$(F_a\times 1)(1\times F_a)(F_b\times 1)=(1\times F_b)(F_a\times 1)(1\times F_a)$ we say that $F$ {\it supports} $R_3(a,a,b)$.  If
$\cal K$ prohibits  $R_3(a,a,b)$ and $F$ satisfies $(F_a\times 1)(1\times F_a)(F_b\times 1)\ne(1\times F_b)(F_a\times 1)(1\times F_a)$
we say that $F$ {\it respects} $R_3(a,a,b)$.

\definition{ If the map $F$ satisifes the following requirements for $\cal K$:
\items
\item{i)}  For each crossing type $C_t$ that is regular in $\cal K$, the component map $F_t$ is fully formed
\item{ii)} $F$ supports all of the allowed dominance relations permitted by $\cal K$ 
\item{iii)} $F$ supports all of the allowed commuting relations permitted by $\cal K$ 
\enditems
\noindent we say that $(X,F)$, or simply $F$ if $X$ is understood, defines a {\it generalised birack} for $\cal K$.}

  If in addition, 
$F$ respects all of the dominance relations that are forbidden by $\cal K$, we say that $F$ is an {\it essential birack} for $\cal K$.

If in addition, for some crossing type $T_a(\Delta)=\Delta$, or equivalently $s_a^y(x) = y \Leftrightarrow s^a_x(y)=x$ then we say that the birack is an {\it $a$-biquandle}, or just a {\it biquandle} if the crossing type is understood.

In most knot theories there is at least one crossing type for which $R_3(a,a,a)$ is permitted.  
Consider a braid-like $R_3$ move describing this dominance and replace each $a$ with $\bar a$.
Since the $R_3$ configuration is symmetrical, the Yang Baxter identities remain valid if we label the
new diagram from the right, rather than the left.  Thus we may 
demonstrate, using the equivalence of $R_3$ moves under the Turaev trick in \cite{F1}, that $a$ also dominates $\bar a$.

Lemma 3.1 shows that we may define biracks in terms of either $S$ or $T$, taking care that the appropriate dominance conditions are satisfied, whereupon $T$ or $S$ respectively is sometimes referred to as the {\it sideways map} for the birack.

Let ${\cal L}(D,X,F)$ be the set of labellings of the diagram $D$ with elements of the label set $X$ that are globally consistent 
at each crossing with the birack $F$.  That is, if an element of ${\cal L}(D,X,F)$ associates the four labels, $x,y,z,w$ with 
a crossing of $D$, then $(z,w) = F(x,y)$ and we say that $F$ satisfies the {\it crossing relation} $(z,w) = (y_x,x^y)$.  If ${\cal L}(D,X,F)$ is non-empty we say that $D$ may be {\it labelled by $(X,F)$}.

As an illustration consider the rack on $\Z_3$ with up operation $x^y=x+1, \hbox{ mod } 3$ and the following labelling of the trefoil diagram. 

\diagram{Labelled  trefoil}

Fig.\ 14 is the closure of a braid.  For any knot theory $\cal K$ there are corresponding braid groups $B_r{\cal K}$, $r \in \N$, with
generators $\alpha_i(a), \alpha_i(\bar a)$ for $i=1, \ldots, r$ and each crossing type $a$ included in 
$\cal K$.  Representative diagrams are braids on $r$ strands in the usual sense, with generalised crossings between strands, whose closure is a braided knot diagram representing an element of $\cal K$. Further details may be found in \cite{BF2}.

If $(X,F)$ is a birack in $\cal K$ then any braid $\beta \in B_r{\cal K}$ defines a bijective function 
$$\beta_\#:X^r\to X^r$$
where $X^r$ the set of column vectors with entries in $X$. Thus, a labelling of the braid closure 
by the birack is defined by a fixed point $\x\in X^r$, where $\beta_\#\x=\x$.

\lemma{If $D'$ is a diagram obtained from $D$ by a permitted $R_2, R_3$ or $R_4$ move, then there is a uniquely defined labelling in ${\cal L}(D',X,F)$ induced by each labelling in ${\cal L}(D,X,F)$.}
{\bf Proof.} Since $(X,F)$ is a birack for $\cal K$, the lemma follows from our observations on the implications of the $R$ moves in sections 4.2.2,
 4.2.3 and 4.2.4.\qed

We call expressions such as $z_{xy_x}, z_{yx^y}, x^{zy^z}, x^{yz_y}, {y_x}^{z_x}, {y^z}_{x^z}$ {\it birack words}; the Kauffman
notation $a_b=a\down{b}$ and $a^b=a\up{b}$ provides the motivation for the use of the term word in this context.

The {\it fundamental birack of $D$, ${\cal F}(D)$, in $\cal K$} is the set of equivalence classes of birack words in generators corresponding to the arcs of $D$ 
under the equivalence relation generated by the crossing relations and the generalised Reidemeister 
$R_2, R_3$ and $R_4$ moves permitted 
by $\cal K$. 

For a given theory $\cal K$, if $L \in {\cal L}(D,X,F)$ there is a birack homomorphism $\phi:{\cal F}(D) \rightarrow (X,F)$; that is
$\phi(x^y) = f_t^{\phi(y)}(\phi(x))$ and $\phi(x_y) = f^t_{\phi(y)}(\phi(x))$ for some crossing type $t$, whose image contains the labels appearing in $L$.


\section{A generalisation of a theorem of Trace}

If two representative diagrams of a classical knot are equivalent under $R_2, R_3$ moves and with no involvement of $R_1$ moves there is a 
{\it regular homotopy} between the diagrams.  A theorem of Trace, \cite{T}, shows that this is dependant on the equality of the {\it writhe and turning number} (defined below). 
In this section we extend this result to the following knot theories that we collectively refer to as 
{\it classical Trace theories}:

\items
\item{i)} Classical knots, ${\cal K}(r)$
\item{ii)} Rotational virtual knots, ${\cal RK}(r,v)$
\item{iii)} The singular knot varieties, ${\cal SK}(r,s)$ and ${\cal SRK}(r,v,s)$ 
\enditems

Rotational virtual knots are studied in detail in \cite{K1}. They are represented as equivalence classes of virtual knot diagrams under the usual classical, virtual and mixed Reidemeister moves with the exception of the virtual $R_1(v)$ move.

These theories share the property that in any sequence of $R$ moves between representative diagrams of the same knot, the only permitted
$R_1$ moves are the classical $R_1(r)$ moves.  

\pagebreak

\subsection{Writhe and turning number}
The {\it writhe}, $w(D)$  of a classical knot diagram $D$, is given by
$$w(D)=\#\{\hbox{ positive crossings}\}-\#\{\hbox{ negative crossings}\}$$
Alternatively, it is the linking number of a parallel longitude defined by the blackboard framing. However it is the first definition which 
can easily be extended to generalised knot diagrams in theories that include classical crossings simply by counting the sum over the 
classical crossings and not involving other crossing types at all.

The writhe is invariant under all $R_2$, $R_3$ and $R_4$ moves involving any crossing types but we can change the writhe to whatever we want by classical $R_1$ moves. In classical knot theory it follows that the writhe is an invariant of regular homotopy.  The $R_4$ 
move in singular knot theories is not an invariant of regular homotopy, and so we make the following more general definition.

\definition{Two diagrams in a generalised knot theory are said to be {\it Trace equivalent} if there is a sequence of $R$ moves connecting
them that involves only $R_2$, $R_3$ and $R_4$ moves.} 

The following lemma follows immediately from the above definition and Lemma 4.2.

\lemma{If $(X,F)$ is a birack for $\cal K$ and $D, D'$ are Trace equivalent, 
then there is a uniquely defined labelling in ${\cal L}(D',X,F)$ induced by each labelling in ${\cal L}(D,X,F)$.\qed}

The {\it turning number} of a generalised knot diagram $D$ is a characteristic of the underlying graph, $|D|$.  We assume that the image 
of the immersed oriented circles that form the unicursal components of $|D|$ are smooth and have a unit tangent at every point. 
The direction of the tangent defines a continuous map from a unicursal component of the $|D|$
to the circle, $S^1$. The resulting integer $n(K)\in\pi_1(S^1)$ corresponding to that map is the {\it turning number} of the component. If the knot has more than one component 
then the (total) turning number is the sum over all components.

We may calculate the turning number of a knot diagram by considering the orientation of the Seifert circles obtained by an orientation
respecting smoothing of the 
crossings of the underlying graph $|D|$.  Attach a sign to each a Seifert circle by assigning $\pm1$ if it is oriented 
clockwise or anti-clockwise, according to a chosen convention of whether an anti-clockwise oriented simple closed curve should have turning number $1$ or $-1$.

\lemma{The turning number of a knot diagram $D$ is equal to the sum of the Seifert circle signs.}
{\bf Proof.} By Seifert smoothing crossings of the flattened virtual knot shadow we create a diagram, $D'$, that is the union of a number of 
disjoint smooth circles, from which we can
determine a turning number based on the map to $S^1$ induced by a unit tangent vector.  The turning number of $D'$ is thus the sum of the 
Seifert circle signs.  If we compare the 
effect on the induced map to $S^1$ of following the Seifert circles rather than the underlying immersion we see that at a crossing, 
on one Seifert circle the unit tangent vector turns to the left(right) by the angle subtended by the adjacent incoming and outgoing arcs.  
On the other Seifert circle the unit tangent vector turns in the opposite direction by the same ammount and so these contributions to the 
turning number cancel out.  The resulting element of $\pi_1(S^1)$ determined by following the components of $D$ or of $D'$ is therefore the same.
\qed

\subsection{Generalised Trace theorem}

\theorem{If $D_1$ and $D_2$  represent the same knot in a classical Trace theory and have the same writhe and turning number then they are Trace equivalent.}
{\bf Proof.}
The diagrams $D_1$ and $D_2$ are related by a sequence of $R$ moves permitted by the theory, based on classical, virtual, and singular 
crossing types as appropriate and {\it not} including the virtual $R_1(v)$ move or the singular $R_1(s)$ move, even where the theory involves these crossing types.  Thus, the only move which change the writhe and turning number are the four shown in Fig.\ 15.
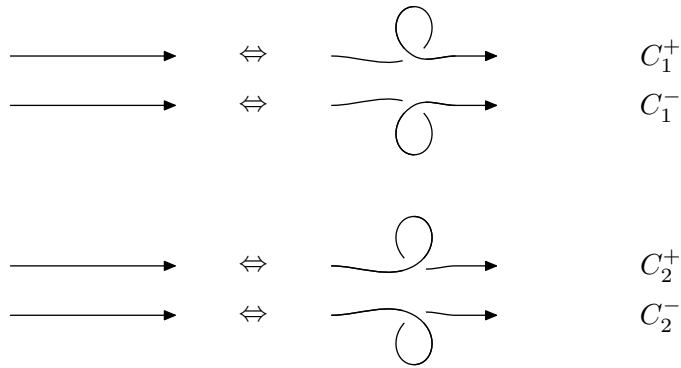
\diagram{Classical $R_1$ moves}
These are paired so that they cancel when next to one another. For example Fig.\ 16 shows $C_1^+$ cancelling with $C_1^-$ using a sequence of $R_2, R_3$ moves, i.e. a Trace equivalence.
\diagram{Cancelling curls}
We now {\it thicken} the curve to make an annulus, see Fig.\ 17.
\diagram{A thickened virtual trefoil}
We now follow the $R$ moves taking $D_1$ to $D_2$. For moves other than classical $R_1$ moves the annulus follows them. When an expansive 
classical $R_1$ curl is encountered then its negative is introduced into the core curve as in Fig.\ 18 where the moves are $C_2^+$ and $C_2^-$.  
\diagram{$C_2^-$ cancelling a $C_2^+$}
On the other hand if a curl is deleted then it will be registered in the annulus such as in Fig.\ 19 which shows the effect of a $C_1^-$ collapse.
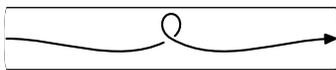
\diagram{The ghost of a $C_1^-$}
We now conclude that there is a Trace equivalence between $D_1$ and $D_2$ with added small curls. Notice that these curls may be moved about
the diagram independently.  A curl may pass a classical or virtual crossing because we have the classical and virtual $R_2$ and $R_3$ moves.
It may pass a singular crossing since we have the dominance moves $R_3(r,r,s)$, $R_3(\bar r,\bar r,s)$ and the move $R_4(r,s)$, as shown for a positive curl in Fig.\ 20.

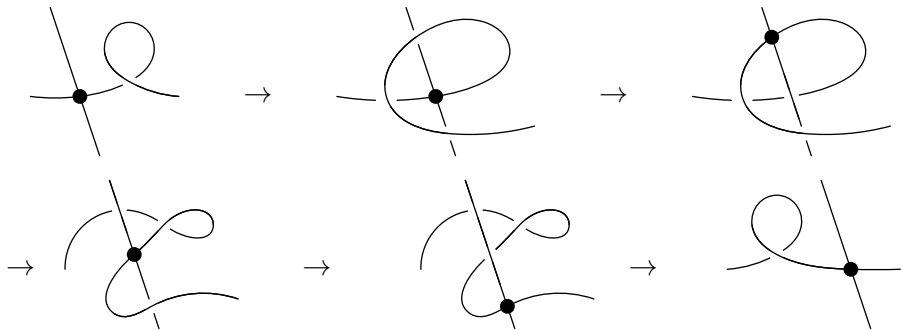
\diagram{Moving past a singular crossing}

Since the writhe and turning number
are preserved by Trace equivalence and are equal for $D_1$ and $D_2$, the curls in $D_2$ may be moved so they cancel in pairs and so disappear.
\qed

We have a similar result for the flat rotational virtual knots, ${\cal RK}(f,v)$ and the corresponding singular variety ${\cal SRK}(f,v,s)$, which collectively we refer to as {\it flat Trace theories}.

\theorem{If $D_1$ and $D_2$  represent the same knot in a flat Trace theory and have the same turning number then they are Trace equivalent.}
{\bf Proof.} The argument follows on the same lines as the previous theorem. Of course the writhe is irrelevant and undefined in this situation. \qed

\subsection{Subdivision and periodic phenomena of ${\cal D}(K)$}

By Theorem 5.5 the graph ${\cal D}(K)$ of all representative diagrams a knot $K$ in one of the classical Trace theories has a two-fold 
subdivision into connected subsets ${\cal D}_{w,n}(K)$, where the vertices of ${\cal D}_{w,n}(K)$ are diagrams with writhe $w$ and turning number $n$ and the edges are Reidemeister moves other than classical $R_1$ moves.

Let $F$ be a birack and let ${\cal L}_{w,n}(K,X,F)$
be the set of labellings of any representing diagram of $K$ by $(X,F)$ with writhe $w$ and turning number $n$. Any diagram in  
${\cal D}_{w,n}(K)$ can be used because there is a bijective equivalence between the labellings of any two diagrams by Lemma 5.3.

There is a lot of redundancy and periodicity in ${\cal D}_{w,n}(K)$ and ${\cal L}_{w,n}(K,X,F)$. The astute reader will have noticed
in Fig.\ 14 that without changing the writhe or the labelling we can change the turning number to $+2$ by looping the bottom arc labelled 
0 above the top arc labelled -1. 

A useful visualisation is to assign a colour to 
${\cal D}_{w,n}(K)$ depending on whether the difference $w-n$ is odd or even.  We may then think of ${\cal D}(K)$ with a chess board colouring
and identify a particular ${\cal D}_{w,n}(K)$ with its corresponding square of a particular colour.  We regard the columns of the chess board 
as having the same writhe and the rows as having the same turning number.

Since the diagrams in ${\cal D}_{w,n}(K)$ are related by $R_2, R_3$ and $R_4$ moves, these moves correspond to edges between diagrams in the
same square and (classical) $R_1$ moves correspond to edges joining diagrams in diagonally adjacent squares of the same colour, since those
moves change both $w$ and $n$ by $1$. Each knot $K$ therefore has representative diagrams in every square of exactly one colour.  

\lemma{For all $w,n$, if $K$ is a knot in the classical Trace theory $\cal K$ and ${\cal D}_{w,n}(K)$ is non-empty, there is a representative diagram in ${\cal D}_{w,n}(K)$ which is the closure of a braid $\beta\in B_r{\cal K}$.}
{\bf Proof.}
The Vogel algorithm \cite{V} for classical knots uses a particular kind of non-parallel $R_2(r)$ moves to convert any knot diagram into the
closure of a classical braid.  In \cite{BF1} it is shown how this algorithm may be extended to any generalised knot theory that includes a regular crossing type. Since only $R_2$ moves are involved in the algorithm it follows that the result is still in ${\cal D}_{w,n}(K)$.  
\qed

Here we regard braids running horizontally from left to right with the first strand at the bottom and the $r$-th strand at the top.  

\lemma{If $D\in {\cal D}_{w,n}(K)$ is the closure of the braid $\beta$, at least one of whose strands closes anti-clockwise (respectively clockwise) 
to form $D$, then there is a diagram $D'\in {\cal D}_{w,n-2}(K)$ (respectively ${\cal D}_{w,n+2}(K)$) that is the closure of the same braid.}
{\bf Proof.}  The Seifert circles in a diagram that is the closure of a braid are in one-to-one correspondence with the braid strands.  
Therefore, if a braid strand is changed from closing anti-clockwise to clockwise, or vice versa, by Lemma 5.4 the turning number changes by 
two, whilst the writhe remains the same.\qed

\corollary{If $D\in {\cal D}_{w,n}(K)$ is the closure of braid $\beta$ then there is a diagram in either ${\cal D}_{w,0}(K)$ or ${\cal D}_{w,1}(K)$) that is the 
closure of the same braid, according to whether $n$ is even or odd.}
{\bf Proof.}  If $n$ is even, repeated application of Lemma 5.4 results in a diagram in ${\cal D}_{w,n+2k}(K)$ for some $k$ that is the closure of 
$\beta$ where all the strands close anticlockwise.  It follows that $\beta$ has an even number, $n+2k$, of strands, so there is a diagram in ${\cal D}_{w,0}(K)$ that is the 
closure of $\beta$.  Similarly, if $n$ is odd.\qed

\theorem{Let $D_1$ and $D_2$ be diagrams of the knot $K$ with the same writhe.  Then for any birack $(X,F)$, 
there is a bijective correspondence between the labellings of $D_1$ and $D_2$ by $(X,F)$.}
{\bf Proof.}  Suppose the turning numbers of $D_1$ and $D_2$ are $n_1$ and $n_2$ respectively and their writhe is $w$; note that $n_1-n_2$ is even.  By Lemma 5.5, there are diagrams
in ${\cal D}_{w,n_1}(K)$ and ${\cal D}_{w,n_2}(K)$ that are the closure of braids $\beta_1$ and $\beta_2$.  By Corollary 5.1 there are 
diagrams $D'_1, D'_2$ in either ${\cal D}_{w,0}(K)$ or ${\cal D}_{w,1}(K)$, depending on whether $n_1, n_2$ are even or odd, that are also
the closure of $\beta_1$ and $\beta_2$.  The diagrams $D'_1, D'_2$ are Trace equivalent and so a labelling of one by $(X,F)$
induces a labelling of the other.  The labellings of $D'_1, D'_2$ are in bijective correspondence with the fixed points of the induced functions 
$\beta_{1\#}:X^{r_1}\to X^{r_1}$ and $\beta_{2\#}:X^{r_2}\to X^{r_2}$, where $r_1,r_2$ 
are the number of strands of the braids $\beta_1, \beta_2$. Thus, there is a bijection between the fixed points of $\beta_{1\#},\beta_{2\#}$
and consequently between the labellings by $(X,F)$ of the diagrams $D_1, D_2$.\qed

In fact, we can go further than Theorem 5.7, as the following lemma shows.

\lemma{Any two diagrams of a knot $K$ in a classical Trace theory with the same writhe have isomorphic fundamental biracks.}
{\bf Proof.}  Any two diagrams $D_1,D_2 \in {\cal D}_{w,n}(K)$ are Trace equivalent, so $D_1,D_2$ are related by $R_2,R_3,R_4$ moves and 
therefore have isomorphic fundamental biracks.  Moreover, two diagrams that are the closure of the same braid have isomorphic fundamental 
biracks, so the result follows from the same argument as in the proof of Theorem 5.7.\qed

Lemma 5.7 means that we can define a stratified fundamental birack ${\cal F}(K)$ for any knot $K$ in a classical Trace theory where, for each writhe $w$ we define ${\cal F}_w(K)$ to be the fundamental birack of any representative diagram of $K$ having writhe $w$.

We end this section by noting that the chess board visualisation provides us with an unexpected result with respect to the writhe and turning number of classical knot diagrams.

\definition{Let $D$ be a classical knot diagram with writhe $w$ and turning number $n$, then the  {\it colour} of $D$ is the odd or even 
parity of the quantity $w-n$.}

\theorem{For any oriented classical knot diagram $D$, the colour of $D$ is odd or even according to whether the number of components of $D$ is odd or even.}
{\bf Proof.} An oriented classical knot diagram may be reduced to a diagram of the unknot using Reidemeister $R_2$ and $R_3$ moves, and changing
positive crossings to negative crossings.  Applying $R_2$ and $R_3$ moves to a diagram in 
${\cal D}_{w,n}(K)$ result in another diagram in ${\cal D}_{w,n}(K)$ so do not change the colour.  Changing the sign of a crossing produces a diagram of a different knot $K'$, with the same turning number but a writhe that is reduced by two, so again, the colour remains the same.
Eventually we arrive at a diagram of the unknot $U$ having the same colour as the original diagram and as all other diagrams of the unknot.  The 
simplest such diagram is a disjoint union of simple loops, which has writhe zero and turning number given by the sum of the turning number of each 
loop.  Thus the colour is the parity of the number of loops counted mod 2.\qed

\section{Finite biracks and the birack polynomial invariant}

In this section, using the previous results, we use finite biracks for labelling knots in a classical Trace theory so we may calculate the 
number of labellings of any given knot.  This leads to an invariant for these knot categories that, since we shall be considering turning 
numbers, we shall consider to be represented by diagrams in the plane.

Let $\phi_w(K,X,F)$ denote the number of labellings in ${\cal L}_{w,n}(K,X,F)$ for a given finite birack $(X,F)$.
Theorem 5.7 tells us that $\phi_w(K,X,F)$ depends only on $w$ and is independent of $n$ and representative diagram.  Referring to our chess 
board analogy for ${\cal D}(K)$, it says that for any knot $K$ and any birack $(X,F)$, the number of labellings by $(X,F)$ of any diagram in a 
square of a given column is the same.

Consider adding a positive classical Markov stabilization move to a braid $\beta$; that is, if $\beta$ is formed of $n$ strands, append the 
term $\sigma_n(r)$ to the braid word.  Take a vector 
$\x\in X^n$ and suppose that the label $y_0$ is associated with the $n$-th strand before the first term of $\beta$.  Let $x_0$ be the label in $\beta_\#\x$ associated with the $n$-th strand after the last term of $\beta$.  
Then in the additional term $\sigma_n(r)$ the birack twitch $T$ takes the pair $(x_0,y_0)$ to $(y_1,x_1)$, say.  If we now add a second stabilization move $\sigma_{n+1}(r)$ to the braid, at this
term the action of the twitch is to take the pair $(x_1,y_1)$ to $(y_2,x_2)$, as shown in Fig.\ 21.
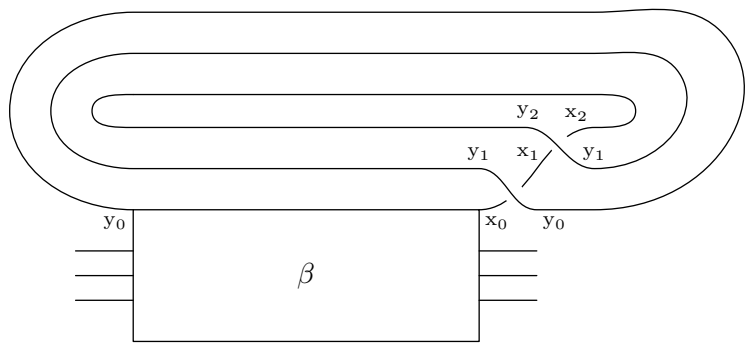
\diagram{Adding stabilization moves to a braid $\beta$}
Continuing in this manner, if $\tau:X^2\rightarrow X^2$ is the bijection $\tau(a,b) = (b,a)$, we see that
the labels after $i$ stabilization moves are given by $W^i(x_0,y_0)$, where $W=\tau\circ T$.  Viewed in terms of the chess board analogy, each successive
Markov $R_1$ move takes us to a diagram in ${\cal D}_{w+1,n+1}(K)$, that is in the adjacent north-east square of the same colour.

A similar argument shows that the labels after $i$ negative Markov stabilization moves at the bottom of the braid, closed clockwise, are given by $W^{-i}$, where $W^{-1}=T^{-1}\circ \tau$.  
Each successive negative stabilization move taking us to a diagram in ${\cal D}_{w-1,n-1}(K)$, that is in the adjacent south-west square of the same colour.

Since $T$ and $\tau$ are bijections, so is $W$ and, since $X$ is finite, $W$ is a permutation and therefore decomposes as a product of cycles.  It 
follows that each pair $(x_0,y_0)$ is part of a cycle whose length, $l$, is such that $W^l$ fixes $(x_0,y_0)$.  We call $l$ the {\it period} of $(x_0,y_0)$ with respect to $T$.

Each fixed point $\x_i$ of $\beta_\#$ yields a pair $x_{i0}=y_{i0}$ in the diagonal of $X^2$.  Each $(x_i,x_i)$ in the diagonal has a period with 
respect to $T$, denoted $l_i$.  
Note that the $l_i$ are not necessarily all the same, since the pairs $(x_i,x_i)$ may lie in cycles of $W$ of different length.

Let $k$ be the least common multiple of the $l_i$ for the diagonal elements of $X^2$, then each $(x_i,x_i)$ is fixed by $W^k$ so a fixed point of $\beta_\#$ induces a fixed point of $\beta\sigma_n\ldots\sigma_{n+k}$.  Thus the number of labellings of the closures of $\beta$
and $(\beta\sigma_n\ldots\sigma_{n+k})_\#$ by the birack $(X,F)$ is the same.

In fact, all we need for $W^k$ to induce a fixed point is for it to preserve the diagonal of $X^2$, it does not need to fix every element of the diagonal.
If some of the diagonal elements are evenly spaced within a cycle factor of $W$, we may use that even spacing as their contribution to the least common
multiple and thereby obtain a minimum value for $k$.

Since $\phi_w(K,X,F)$ is well defined the power series
$$\sum_{-\infty}^{+\infty} \phi_w(K,X,F)t^w$$
is well defined.  We have seen that there is a repeating pattern of $k$ coefficients to this power series, where $k$ is the minimum integer such that $W^k$ fixes the diagonal of $X^2$, so the finite polynomial
$$f(K,X,F)=\sum_{0}^{k-1} \phi_w(K,X,F)t^w$$
is therefore a knot invariant, called the {\it birack polynomial} of $K$ with respect to $(X,F)$.

Each element of ${\cal L}(D,X,F)$ may be viewed as the labelling of $D$ by a homomorphism from the fundamental birack ${\cal F}(D)$ to $(X,F)$
and we have seen in Lemma 5.7 that any diagram of $K$ with writhe $w$  has a fundamental birack isomorphic to ${\cal F}_w(K)$.  We may refine
the birack polynomial by writing $\phi_w(K,X,F) = |\hbox{Hom}({\cal F}_w(K),(X,F))|$ and noting the size of the image of each homomorphism, 
provided this refinement is well defined.  
We note, in particular, that in some cases, not all of the elements of a homomorphism image are required to label $D$.

However, the image of a birack homomorphism $\phi:(X,F)\rightarrow(X',F')$ is a sub-birack of the co-domain, since it inherits the properties required
by definition 4.1 from the co-birack $(X',F')$ and is closed since $\phi(x^y) = f_t^{\phi(y)}(\phi(x))$ and $\phi(x_y) = f^t_{\phi(y)}(\phi(x))$
for each component map $F_f$ of $F$.   Moreover we have the following lemma.

\lemma{For each labelling $L \in {\cal L}(D,X,F)$, there is a unique smallest sub-birack of $(X,F)$ containing the labels that appear in $L$.}
{\bf  Proof.}  We may describe the up and down actions of the switch associated with a finite birack as a pair of matrices $U,D$ where the $(x,y)$ entries
are $y^x$ and $y_x$ respectively, meaning the rows of $U,D$ are pemutations.

A subset $A\subset X$ forms a sub-birack $(A,S)$ of $(X,S)$ if, and only if, for both $U$ and $D$, each row corresponding to an 
element of $A$ contains an element of $A$ in each column corresponding to an element of $A$.  It follows that, if $A,B,C$ are mutually disjoint subsets of 
$X$ and both $(A\cup B,S)$ and $(A \cup C,S)$ are sub-biracks of $(X,S)$, then $(A,S)$ is a sub-birack of $(X,S)$.  This is because
no row of $U,D$ corresponding to an element of $A$ can have an element of $A$ in a column corresponding to an element of $B$ or $C$, since $(A \cup C,S), 
(A\cup B,S)$ respectively are biracks.\qed

Any two diagrams in ${\cal D}_{w,n}(K)$ are Trace equivalent so by Lemma 5.3, if $D,D' \in {\cal D}_{w,n}(K)$ and the smallest image of a birack 
homomorphism containing the labels of $L \in {\cal L}(D,X,F)$  is the sub-birack $(A,F)$, then  $L$ induces a unique labelling, $L'$, of $D'$ by $(A,F)$.  
Moreover the labels in $L'$ cannot lie in a proper sub-birack $(B,F)$ of $(A,F)$, since the only way that fewer labels may appear in $L'$ than in 
$L$ is by nature of a negative $R_2$ move in the sequence of $R$ moves connecting $D$ and $D'$.  The labels attached to the edges of a 
Reidemeister 2 bigon necessarily belong to the same birack as the other labels attached to the bigon crossings, since they are the result 
of the up and down action between those labels.

Thus, for $i=1,\ldots,r$,  for any two diagrams in ${\cal D}_{w,n}(K)$ the number of labellings whose smallest containing sub-birack is 
the image of a homomorphism from ${\cal F}_w(K)$ of size $i$ is the same.  We may therefore write 
$|\hbox{Hom}({\cal F}_w(K),(X,F))| = \sum_{i=1}^r n_i$, where $|X|=r$ and there are $n_i$ homomorphisms in $\hbox{Hom}({\cal F}_w(K),(X,F))$
whose image involves $i$ elements of $X$ and is the smallest sub-birack of $(X,F)$ containing all the labels in an element of ${\cal L}(D,X,F)$.

Now we proceed as in the proof of Theorem 5.7 to show that the partition $\{n_1,\ldots,n_r\}$ of $|\hbox{Hom}({\cal F}_w(K),(X,F))|$ is well 
defined when determined by any diagram $D$ with writhe $w$.  Thus we can write $\phi_w(K,X,F) = \sum_{i=1}^r n^w_i$ without ambiguity and 
can write the birack polynomial as a two variable polynomial.

$$f(K,X,F)=\sum_{w=0}^{k-1} \big(\sum_{i=1}^r n^w_i s^{i-1}\big)t^w $$

\subsection{Examples}

In the lists of finite biracks given at www.layer8.co.uk/maths/biquandles we have the following:

$R^5_{40}$ from racks of size 5 that are not quandles: $$U=((1 3) , (4 5) , (1 3) , (2 5) , (2 4))\quad D=\iota$$ \break
$R^6_{114}$ from racks of size 6 that are not quandles: $$U=((2 3)(4 5 6) , (2 3) , (2 3) , (1 6 5)(2 3) , (1 4 6)(2 3) , (1 5 4)(2 3))\quad D=\iota$$ \break
$BR^6_{125}$ from quandle-related biracks of size 6 that are not biquandles and not racks:
$$U=((3 4 6) , (3 4 6) , (1 5 2) , (1 5 2) , (3 4 6) , (1 5 2))$$ $$D=((1 2 5)(3 4 6) , (1 2 5)(3 4 6) , (4 6) , (3 6) , (1 2 5)(3 4 6) , (3 4))$$

$$\vbox{ \offinterlineskip \halign{\strut
\vrule \hfil \quad # \quad \hfil \vrule & \hfil\ #\ \hfil \vrule & \hfil\ #\ \hfil \vrule  & \hfil\ #\ \hfil \vrule  
& \hfil\ #\ \hfil \vrule & \hfil\ #\ \hfil \vrule  & \hfil\ #\ \hfil \vrule \cr 
\noalign{\hrule} 
& & \multispan3 \hfil Birack polynomial\hfil \vrule \cr
\omit\vrule & \omit\vrule & \omit\hrulefill\vrule & \omit\hrulefill\vrule & \omit\hrulefill \cr
\omit\vrule height2pt \hfil \vrule height2pt&                 &     &     & \cr
 knot & braid word         & $R^5_{40}$  & $R^6_{114}$ & $BR^6_{125}$ \cr 
\noalign {\hrule} 
  unknot   &               &  $3t + 5$      &  $4t+6$        &  $3t^2 + 3t +6$     \cr
 3.1 & $\s_1\s_1\s_1 $              &  $9t + 11$   &  $16t+ 18$  &  $9t^2 +9t +12$   \cr 
 4.1 & $\s_1\s^{-1}_2\s_1\s^{-1}_2 $         &  $3t + 5$     & $16t+18$  &  $3t^2 +3t +6$     \cr
\noalign{\hrule}}}$$

\centerline{
$$\vbox{ \offinterlineskip \halign{\strut
\vrule \hfil \quad # \quad \hfil \vrule & \hfil\ #\ \hfil \vrule & \hfil\ #\ \hfil \vrule  & \hfil\ #\ \hfil \vrule  
& \hfil\ #\ \hfil \vrule & \hfil\ #\ \hfil \vrule  & \hfil\ #\ \hfil \vrule \cr 
\noalign{\hrule} 
& & \multispan3 \hfil Refined birack polynomial\hfil \vrule \cr
\omit\vrule & \omit\vrule & \omit\hrulefill\vrule & \omit\hrulefill\vrule & \omit\hrulefill \cr
\omit\vrule height2pt \hfil \vrule height2pt&                 &     &     & \cr
 knot & braid word         & $R^5_{40}$  & $R^6_{114}$ & $BR^6_{125}$ \cr 
\noalign {\hrule} 
  unknot   &               &  \smaths{3t+(2s+3)}      &  \smaths{4t+(2s+4)}        &  \smaths{3t^2+3t+(3s^2+3)}     \cr
 3.1 & $\s_1\s_1\s_1 $              &  \smaths{(6s^2+3)t+(6s^2+2s+3)}   &  \smaths{(12s^3+4)t+(12s^3+2s+4)}  &  \smaths{(6s^2+3)t^2+(6s^2+3)t+(9s^2+3)}   \cr 
 4.1 & $\s_1\s^{-1}_2\s_1\s^{-1}_2 $         &  \smaths{3t+(2s+3)}     & \smaths{(12s^3+4)t+(12s^3+2s+4)}  &  \smaths{3t^2+3t+(3s^2+3)}     \cr
\noalign{\hrule}}}$$
}

In \cite {B}, Bigelow presents two nontrivial braids, $\beta_1, \beta_2$, that lie in the kernel of the Burau representation of braids for 
braids on five and six strands respectively.  The two braids are defined as follows
$$\beta_1 = [\psi_1^{-1}\sigma_4\psi_1,\psi_2^{-1}\sigma_4\sigma_3\sigma_2\sigma_1^2\sigma_2\sigma_3\sigma_4\psi_2]$$
\noindent where $\psi_1=\sigma_3^{-1}\sigma_2\sigma_1^2\sigma_2\sigma_4^3\sigma_3\sigma_2$  and 
$\psi_2 =\sigma_4^{-1}\sigma_3\sigma_2\sigma_1^{-2}\sigma_2\sigma_1^2\sigma_2^2\sigma_1\sigma_4^5$
$$\beta_2 = [\psi_1^{-1}\sigma_3\psi_1,\psi_2^{-1}\sigma_3\psi_2]$$
\noindent where $\psi_1=\sigma_4\sigma_5^{-1}\sigma_2^{-1}\sigma_1$ and $\psi_2 =\sigma_4^{-1}\sigma_5^2\sigma_2\sigma_1^{-2}$
and where $[\psi_1,\psi_2]= \psi_1^{-1}\psi_2^{-1}\psi_1\psi_2$.

The closure of these braids, $\tilde\beta_1, \tilde\beta_2$, have five and six components repsectively and, as shown in the following table, are distinguished from the unlink on five and 
six components by the birack polynomial.

$$\vbox{ \offinterlineskip \halign{\strut
\vrule \hfil \quad # \quad \hfil \vrule & \hfil \quad # \quad \hfil \vrule & \hfil\ #\ \hfil \vrule\cr 
\noalign{\hrule} 
& \hfil Birack polynomial\hfil & \hfil Refined birack polynomial\hfil\cr
 link          & $R^5_{40}$ & $R^5_{40}$ \cr 
\noalign {\hrule} 
  5-component unlink   &  \smaths{1875t+3125} & \smaths{(1392s^4+480s^2+3)t+(2220s^4+870s^2+32s+3)} \cr
 $\tilde\beta_1$ & \smaths{1443t+2549} & \smaths{(960s^4+480s^2+3)t+(1644s^4+870s^2+32s+3)} \cr
  6-component unlink   &  \smaths{9375t+15625} & \smaths{(7920s^4+1452s^2+3)t+(12840s^4+2718s^2+64s+3)} \cr
 $\tilde\beta_2$ & \smaths{3567t+7273} & \smaths{(2112s^4+1452s^2+3)t+(4488s^4+2718s^2+64s+3)} \cr
\noalign {\hrule} 
}}$$

\section{Further work}

We have shown that, for the classical Trace knot theories, the fundamental birack of any two diagrams of a knot $K$
with the same writhe are isomorphic, which gives us a stratified fundamental birack for the knot itself.  
We have constructed one invariant from this birack and it would be interesting to understand how
it relates to other invariants. Its relationship to the Jones polynomial and quantum link invariants that 
are, for the non-singular theories, invariants of regular homotopy are of particular interest, as are the Vassiliev invariants
for the singular theories.  We also have a Trace theorem for flat virtual and singular flat virtual knots but without a
stratification.  Nevertheless, it seems likely that there are further insights into these categories to be found via the 
chess board visualisation of ${\cal D}(K)$.

\section{Bibliography}

[B] S. Bigelow.  The Burau representation is not faithful for n=5.  Geometry and Topology Vol 3, 1999, pp 397-404.

[BF1] A. Bartholomew, R. Fenn. Alexander and Markov theorems for generalized knots, I. Journal of Knot Theory and Its Ramifications 
Vol. 31, No. 08, 2022.

[BF2] A. Bartholomew, R. Fenn. Alexander and Markov theorems for generalized knots, II. Journal of Knot Theory and Its Ramifications 
Vol. 31, No. 08, 2022.

[F1] R. Fenn. Generalised biquandles for generalised knot theories.  New Ideas in Low Dimensional Topology. March 2015, 79-103

[F2] R. Fenn, Combinatorial Knot Theory, Series on Knots and Everything volume 76, World Scientific Books, 2024.

[FJK] R. Fenn, M. Jordan-Santana and L. Kauffman.  Biquandles and virtual links. in Topology Appli. 145 (2004) 157–175.

[HN] A.  Henrich and S.  Nelson.  Semiquandles and flat virtual knots. Pacific J. Math. {248} (2010), no. 1, 155--170. 

[K1] L. Kauffman. Rotational virtual knots and quantum link invariants.  Journal of Knot Theory and Its Ramifications, Vol. 24, No. 13, 1541008 (2015)

[K2] L. Kauffman. Multi-virtual knot theory.  arXiv:2409.07499v1

[T] B. Trace.  On the Reidemeister moves of a clasical knot. Proc. Am. Math. Soc. Vol 89, No. 4, 1989.

[V] P. Vogel. Representation of links by braids, a new algorithm.  Comment. Math. Helvetici 65 (1990), 104-113.

\bye